\documentclass[12pt]{amsart}
\usepackage[centertags]{amsmath}
\usepackage{graphicx}
\usepackage{amsfonts,amssymb,latexsym,epsfig,amsthm}
\usepackage{fullpage}

\newtheorem{theorem}{Theorem}

\newtheorem{claim}[theorem]{Claim}

\newtheorem{corollary}[theorem]{Corollary}

\newtheorem{lemma}[theorem]{Lemma}

\newtheorem{proposition}[theorem]{Proposition}
\newtheorem{remark}[theorem]{Remark}

\input{tcilatex}
\newcommand{\be}{\begin{equation}}
\newcommand{\ee}{\end{equation}}
\newcommand{\ba}{\begin{eqnarray}}
\newcommand{\ea}{\end{eqnarray}}
\newcommand{\ban}{\begin{eqnarray*}}
\newcommand{\ean}{\end{eqnarray*}}
\renewcommand{\o}[2]{\frac{#1}{#2}}

\newcommand{\intl}{\int\kern-9pt\hbox{$\backslash$}}

\begin{document}
\title{On a Conformal Gauss-Bonnet-Chern Inequality for LCF manifolds and
related topics}
\author{Hao Fang}
\address{Courant Institute of Mathematical Sciences, New York University}
\email{haofang@cims.nyu.edu}
\date{February 20, 2004}

\begin{abstract}
In this paper, we prove the following two results:

First, we study a class of conformally invariant operators $P$ and their
related conformally invariant curvatures $Q$ on even-dimensional Riemannian
manifolds. When the manifold is locally conformally flat(LCF) and compact
without boundary, $Q$-curvature is naturally related to the integrand in the
classical Gauss-Bonnet-Chern formula, i.e., the Pfaffian curvature. For a
class of even-dimensional complete LCF manifolds with integrable $Q$%
-curvature, we establish a Gauss-Bonnet-Chern inequality.

Second, a finiteness theorem for certain classes of complete LCF four-fold
with integrable Pfaffian curvature is also proven. This is an extension of
the classical results of Cohn-Vossen and Huber in dimension two. It also can
be viewed as a fully non-linear analogue of results of Chang-Qing-Yang in
dimension four.
\end{abstract}

\maketitle

\section{Introduction}

\setcounter{equation}{0} \setcounter{theorem}{0}

\label{1}

Let $M$ be a Riemannian manifold of even dimension $n$, with a Riemannian
metric $g$. Denote $[g]=\{e^{2w}g;\ w\in C^{\infty }(M)\}$ as the conformal
metric class determined by $g$. It is known (Cf.~\cite{B} \cite{FG2}) that
there exist local curvatures $Q_{2k}$ (with $2k\leq \dim M$) which satisfy
certain conformal transformation laws if the Riemannian metric varies in the
conformal metric class $[g]$. The most interesting one is $Q_{n}=Q_{\dim M}$%
. Denoted also as $Q$ for future convenience, it satisfies the following
transformation law:
\begin{equation}
Q_{w}\ dv_{g_{w}}=(Q_{0}+P_{n}w)\ dv_{g_{0}},  \label{1.1}
\end{equation}%
where $g_{0}$ and $g_{w}=e^{2w}g_{0}$ are both in $[g]$ and $P_{n}$ is an $n$%
-th order linear elliptic operator. In recent years, significant progress
has been made in the study of $P_{n}$ and $Q_{n}$, for $n=2$ and $n=4$. It
is known to be closely related to the theory of partial differential
operators and spectral invariants. For more details and background, see
Section~\ref{2}.

For a Riemannian manifold $(M,g)$ of dimension $2m,$ an important
characteristic class, Pfaffian invariant, is defined as the $m^{\text{th}}$
Chern class, $c_{m}(M),$ hence it can be represented as a curvature
invariant by the standard Chern-Weil theory.

If we assume that $M$ is locally conformally flat(LCF), it is an interesting
fact that $Q_{n}$ is a multiple of the Pfaffian of the metric modulo a
divergence term. This can be proved by applying a result of
Branson-Gilkey-Pohjanpelto in invariant theory(Cf.~\cite{BGP}). Thus, if $M$
is compact without boundary, the Gauss-Bonnet-Chern theorem gives:
\begin{equation}
C_{n}\int\limits_{M}{Q}_{g}{dv_{g}}=\mathrm{Euler}(M),  \label{GBC}
\end{equation}%
with $C_{n}={\frac{{1}}{{((n-2)!!)^{2}|S^{n-1}|}}}$. (Notice that $%
(2k)!!=(2k)(2k-2)\cdots 2$ for a positive integer $k$ and $0!!=1;$ $|S^{n-1}|
$ denotes the volume of the standard $(n-1)$-sphere of radius $1$.) Here $%
\mathrm{Euler}({M)}$ denotes the Euler Characteristic of $M$, which is a
topological invariant of the manifold.

One of the goals of this paper is to extend the above-mentioned formula to
certain complete LCF manifolds. There are several known results in low
dimensional cases. In dimension two, a classical result for complete open
surfaces by Cohn-Vossen~\cite{CV} and Huber~\cite{H} shows a
Gauss-Bonnet-Chern-type inequality is valid for complete surfaces with
integrable $Q_2$ (which is exactly the Gaussian curvature in dimension two).
In dimension four, Chang, Qing and Yang~\cite{CQY1,CQY2} extended this
inequality to certain complete LCF manifolds with integrable $Q_4$. In this
paper, the general even dimensional case will be considered and the
following will be proven:

\begin{theorem}
\label{01} Assume $(M,g)$ is a complete LCF manifold with finitely many
conformally flat ends with
\begin{equation*}
\int_{M}|Q_{g}|dv<\infty .
\end{equation*}%
If, near the ends, the scalar curvature $R_{g}$ satisfies
\begin{equation*}
R_{g}\geq 0,
\end{equation*}%
then
\begin{equation*}
C_{n}\int_{M}Q_{g}\ dv\leq \mathrm{{Euler}(M),}
\end{equation*}%
where $C_{n}={\frac{{1}}{{((n-2)!!)^{2}|S^{n-1}|}}}$.
\end{theorem}

With the exception of the work of Cheeger-Gromov~\cite{CG} on manifolds with
bounded geometry, and the work of Greene-Wu~\cite{GW} on the complete
4-folds with positive sectional curvature, there is little known about
extensions of the original Gauss-Bonnet-Chern formula for complete manifolds
in higher dimensions. Theorem~\ref{01} suggests that $Q_n$, obtained by
adding a divergence term to the Pfaffian curvature, should be the right
integrand to consider for complete LCF manifolds.

Theorem~\ref{01} is proved by first analyzing the model problem where $M$ is
$\mathbb{R}^{n}$. For the model problem, a geometric averaging argument
further reduces the metrics to rotationally symmetric metrics on $\mathbb{R}%
^{n}$, for which a uniqueness result of the conformal factor is proved by
solving the ODE induced from (\ref{1.1}). The general case is then derived
from the model problem by a gluing argument.

Theorem~\ref{01} can be applied to study the conformal compactification of
certain LCF manifolds, as in~\cite{H} and~\cite{CQY2}. See~\cite{F} for some
details. This will be addressed in a separate paper.

During the course of proving Theorem~\ref{01}, we closely study the
conformal transformation law (\ref{1.1}) of the $Q$ curvature, which is a
linear elliptic PDE with respect to the background metric. Interestingly,
some of the techniques we employ are also effective for various non-linear
problems. In particular, we study the Pfaffian curvature of a complete LCF
four-fold, which satisfies a fully non-linear conformal transformation law
of Monge-Ampere type. (See Section~\ref{5} for more details.)

Hence, the second part of this paper is a generalization of the main result
of~\cite{H, CQY2} in another direction; namely, we consider the
compactification of LCF manifolds.Under some local curvature conditions, we
prove the following finiteness result in dimension four:

\begin{theorem}
\label{03} Let $\Omega $ be a four-fold with a LCF metric. If there exist
constants $C$ and $C^{\prime }$ such that
\begin{equation*}
C\geq R_{g}\geq C^{\prime }>0,\ \Vert \nabla _{g}R_{g}\Vert _{g}\leq C,\
\mathrm{Ric}_{g}\geq -Cg
\end{equation*}%
and
\begin{equation*}
\int\limits_{\Omega }{|\mathrm{Pfaff}_{g}|dv_{g}}<\infty ,
\end{equation*}%
then $\Omega =S^{4}\backslash \{p_{1},\cdots ,p_{k}\}$ for some $p_{i}\in
S^{4}$ ($i=1,\cdots ,k$).
\end{theorem}

In proving Theorem~\ref{03}, the local divergence structure of the Pfaffian
plays an important role, overcoming the difficulties caused by the lack of
linear transformation laws of the Pfaffian.

We would like to comment that there have been extensive studies on the
geometric significance of the corresponding fully non-linear equation in
general metric situation. Especially in dimension four, Chang-Gursky-Yang
proved a conformal sphere theorem~\cite{CGY2}.\ See also \cite{CY6} for
references on further developments.

In a seperate paper, we would like to address the general dimensional cases
of this compactness problem.

This paper is organized as follows: in Section~\ref{2}, some preliminary
facts about conformally invariant operators and curvatures are given. In
Section 3, we prove Theorem~\ref{01}. In Section~\ref{5}, we prove Theorem~%
\ref{03}.

\emph{Acknowledgment}. This material represents part of the author's
doctoral dissertation at Princeton University, 2001. The author would like
to thank his thesis advisor, Alice Chang, for support and guidance. He also
wishes to thank Paul Yang, Jeff Viaclovsky and many others for their
interest in the work and for their helpful discussions.

\section{Conformally invariant operators and curvatures}

\setcounter{equation}{0} \setcounter{theorem}{0} \label{2}

In this Section we give a reivew on the conformally invariant operators and
curvatures.

\subsection{The general metric case}

Let $M$ represent a Riemannian n-fold with a fixed Riemannian metric $g_{0}$%
. Any metric $g$ in the conformal class $[g_{0}]$ can be expressed as $%
g=g_{w}=e^{2w}g_{0}$, where $w$ is a smooth function on $M$. It is therefore
true that the metrics in $[g_{0}]$ can be endowed with an affine structure
modelled after $C^{\infty }(M)$, the linear space of smooth functions on $M$%
. Let $R_{g}$ and $\mathrm{Ric}_{g}$ be the scalar curvature and the Ricci
curvature of $g$, respectively. For future convenience, we define a
symmetric quadratic form,\ the Schouten tensor:
\begin{equation}
A_{g}={\frac{1}{{n-2}}}(\mathrm{Ric}_{g}-{\frac{R_{g}}{{2(n-1)}}}g).
\label{A_g}
\end{equation}%
Denote $\sigma _{k}=\sigma _{k}(A_{g})$ to be the $i$-th symmetric
polynomial of the eigenvalues of $A_{g}$. That is, if we denote $\lambda
_{1},\cdots ,\lambda _{n}$ as eigenvalues of $A_{g}$,
\begin{equation}
\sigma _{k}=\sigma _{k}(A_{g})=\sum\limits_{1\leq i_{1}<\cdots <i_{k}\leq
n}(\prod\limits_{i=1}^{k}\lambda _{i_{j}}).  \label{sigma}
\end{equation}

In particular, we have $\sigma_1={\frac{R_g}{{2(n-1)}}}$, which is also
denoted as $J$ for future convenience.

A conformally invariant operator, $P_{2k}$, is a $2k$-th order partial
differential operator acting on $C^{\infty }(M)$, such that under a
conformal change of the metric $g_{w}=e^{2w}g_{0}$, it obeys the following
transformation law:%
\begin{equation}
P_{2k,w}f=e^{-\alpha w}P_{2k,0}(e^{\beta w}f),  \label{P_{2k}}
\end{equation}%
for some real $\alpha $ and $\beta $. Note that we use subscripts to
indicate the metrics used.

Using the construction of Fefferman-Graham~\cite{FG},
Graham-Jenne-Mason-Sparling~\cite{GJMS} showed that if $\dim M=n$ is odd, $%
P_{2k}$ exists for $(\alpha ,\beta )=(\o {1}{2}{n}+k,\o {1}{2}{n}-k)$, with $%
k$ being any positive integer. If $n$ is even, $P_{2k}$ exists for $(\alpha
,\beta )=(\o {1}{2}{n}-k,\o {1}{2}{n}+k)$, with $k\leq \o {1}{2}{n}$. In
either case, $P_{2k}$ can be assigned the same symbol as that of $\Delta
^{k} $. Furthermore, if $g$ is locally Euclidean,
\begin{equation}
P_{2k}=\Delta ^{k}.  \label{P-flat}
\end{equation}

Except for low order cases, the general expression of $P_{2k}$ is unknown.
Explicit formulae of $P_{2k}$ on $S^{n}$ has appeared in Branson~\cite{B1};
see also Beckner~\cite{Be}. For inductive expressions of the conformally
invariant operators, see~\cite{GJMS, W1, W2}.

Recently, Alexakis has given some general description of the conformal
invariant operators, see~\cite{A}.

If $(\alpha, \beta)\neq (\dim M,0)$, we define
\begin{equation}
Q_{2k}\equiv P_{2k} 1.  \label{Q_2}
\end{equation}
It is clear that $Q_{2k}$ depends only locally on the Riemannian metric.
Hence, by invariance theory, it is fully determined by the curvature tensor
and its covariant derivatives.

For example, if $n>2$, $P_{2}$ is the well-known conformal Laplacian; $Q_{2}$
is then a multiple of scalar curvature. The famous Yamabe problem studies
the existence of constant scalar curvature metric in any conformal metric
class on a compact closed manifold. It was settled by Yamabe, Trudinger,
Aubin and Schoen by using techniques of calculus of variations and studying
the PDE induced from (\ref{Q_2}) (for $k=1$). See~\cite{LP} for more details
and complete references.

However, the most interesting case occurs when $\dim M=n=2m$ is even and $%
(\alpha, \beta)= (n,0)$. From physical considerations (Cf.~\cite{De}, for
example), it is natural to ask if there exists a local curvature invariant $%
Q_n$ satisfying the following conformal transformation law:
\begin{equation}
Q_{n,w}= e^{-nw}(Q_{n,0}+P_{n,0} w).  \label{Q_n}
\end{equation}

In dimension two, it is easy to see that
\begin{eqnarray*}
P_{2} &=&\Delta _{g}, \\
Q_{2} &=&\o {1}{2}R
\end{eqnarray*}%
satisfy (\ref{P_{2k}}) and (\ref{Q_n}). There have been extensive studies
for the geometry of the scalar curvature of closed Riemann surfaces. The
Nirenberg problem asks which functions on a Riemann surface can be
prescribed as the scalar curvature (\textit{i.e.} $Q_{2}$) of a metric in a
given conformal class. On the other hand, from the view of calculus of
variations, $P_{2}$ and $Q_{2}$ are closely related to extremals of the zeta
functional determinant of the Laplacian (\textit{i.e.} $P_{2}$). With the
delicate analytic tools developed, many deep geometric results have been
obtained by Onofri~\cite{O}, Trudinger~\cite{T}, Moser~\cite{M},
Kazdan-Warner~\cite{KW}, Chang-Yang~\cite{CY1,CY2}, Chang-Liu~\cite{CL},
Osgood-Phillips-Sarnak~\cite{OPS1,OPS2} and others. See also~\cite{C} for a
survey.

In dimension four, Paneitz~\cite{P} proved that the following $P_{4}$ and $%
Q_{4}$ satisfy (\ref{P_{2k}}) and (\ref{Q_n}):
\begin{equation}
P_{4}=\Delta _{g}^{2}+\delta (\o {2}{3}R_{g}g-2\mathrm{Ric}_{g})d,
\label{pan4}
\end{equation}%
\begin{equation}
Q_{4}=\o {1}{6}(-3\Vert \mathrm{Ric}\Vert _{g}^{2}+R_{g}^{2}-\Delta
_{g}R_{g}),  \label{P2}
\end{equation}

with $\delta $ being the adjoint operator of $d$ with respect to $g$. In
analogy to the two-dimensional case, it is interesting to study the problem
of prescribing $Q_{4}$ curvature for a given four-manifold as well as the
properties of $Q_{4}$, a 4-th order linear elliptic operator. Furthermore, $%
P_{4}$ and $Q_{4}$ naturally appear in the variation of the functional
determinants of certain conformally invariant operators. Extensive studies
on the analysis and the geometry of the $P_{4}$ and $Q_{4}$ have been
carried out by Beckner~\cite{Be}, Branson-Chang-Yang~\cite{BrCY}, Chang-Yang~%
\cite{CY5}, Chang-Gursky-Yang~\cite{CGY}, Gursky~\cite{G} and many others.
See also~\cite{CY3} and~\cite{CY6} for surveys.

In~\cite{B}, Branson proved the existence of $Q_{n}$ curvature satisfying (%
\ref{Q_n}) for arbitrary even dimensions. Recently,\ Graham and Zworski~\cite%
{GZ} has given a different proof, which was later\ greatly simplified in
\cite{FG2}. The new approach, which is partly based on the fundamental work
of Fefferman and Graham~\cite{FG} on the construction of ambient metric, has
also inspired many related works~\cite{FH}.

However, due to the complicated nature of the expression of $Q_{n}$ (as well
as that of $P_{n}$) for $n$ large, except for a discussion on $Q_{n}$ for
metrics in the standard conformal metric class of $S^{n}$~\cite{CY4}, few
results have been obtained. See~\cite{GP} for dimension eight computation by
using the tractor calculus technique.

Notice that the pair $(P_{n},Q_{n})$ is not unique for $n\geq 4$ in general.
For example, denote $W_{g}$ as the Weyl tensor of the metric. Given a pair $%
(P_{n},Q_{n})$ satisfying the above-mentioned relationship, it is easy to
check that $(P_{n}+c_{1}\Vert W\Vert ^{\o {1}{2}n},Q_{n}+c_{2}\Vert W\Vert ^{%
\o {1}{2}n})$ satisfies the same relations for arbitrary real $c_{1}$ and $%
c_{2}$. See~\cite{A} for a structure theorem of general conformally
invariant curvatures.

\subsection{The LCF metric case}

We now restrict to the case where the metric is locally conformally flat
(LCF). That means, in local coordinates, the metric can be represented as $%
g=e^{2w}g_{0}$, where $g_{0}$ is the standard Euclidean metric. Since the
curvature tensor of the flat metric $g_{0}$ vanishes, one gets $Q_{n,0}=0$.
Thus, by (\ref{P-flat}) and (\ref{Q_n}), we have
\begin{eqnarray*}
P_{n,w} &=&e^{-nw}\Delta ^{m}, \\
Q_{n,w} &=&e^{-nw}\Delta ^{m}w,
\end{eqnarray*}%
where $\Delta $ is the Laplace operator with respect to the flat metric. $%
Q_{n}$ is thus uniquely determined when the metric is LCF. However, as
mentioned in the previous Subsection, the explicit expressions of $P_{n}$
and $Q_{n}$ using the Riemann curvature tensor and its covariant derivatives
are difficult to obtain for higher dimensional cases.

We fix local coordinates $\{x_1,\cdots ,x_n\}$ such that $%
g_{0,ij}=\delta_{ij}$. Hence, locally $g_{ij}= e^{2w}\delta_{ij}$. Under
this coordinate system, the following well-known formulae hold:

\begin{equation}
R_g=-2(n-1)e^{-2w}(\Delta w+{\frac{{n-2}}{2}}\|\nabla w\|^2);
\label{scalar-curvature}
\end{equation}
\begin{equation}
\mathrm{Ric}_{g,ij}=(2-n)w_{ij}-\Delta w \delta_{ij}+(n-2)(w_iw_j-\|\nabla
w\|^2\delta_{ij});  \label{ricci}
\end{equation}
\begin{equation}
A_{g,ij}=-w_{ij}+w_iw_j-\o {1}{2} \|\nabla w\|^2 \delta_{ij},  \label{aij}
\end{equation}
where $w_i={\frac{\partial }{\partial{x_i}}}w=\nabla_i w$ is the derivative
with respect to the flat metric.

We also denote $\mathrm{Pfaff}=c_{m}(M,g)$ as the Pfaffian of the metric $g$%
, with the normalization so that for a closed manifold $M$, the
Gauss-Bonnet-Chern theorem reads:
\begin{equation}
\int\limits_{M}{\mathrm{Pfaff}}_{g}{\ dv_{g}}=\mathrm{Euler}(M).  \label{gbc}
\end{equation}%
Pfaffian invariants for LCF metrics can be expressed as a contraction of the
Schouten tensor as follows (Cf. [V]):
\begin{equation}
\text{Pfaff}_{g}={\frac{{1}}{{((n-2)!!)^{2}|S^{n-1}|}}}\sigma _{m}(A_{g}).
\label{pfaff}
\end{equation}%
Notice that the right hand side of (\ref{pfaff}), when viewed as an
expression of the conformal factor $w$ by (\ref{aij}), is fully non-linear.

The following result is a consequnece of a theorem of Branson, Gilkey and
Pohjanpelto~\cite{BGP}:

\begin{proposition}
\label{11} If $(M,g)$ is an LCF manifold, then

\begin{equation}
{\frac{{1}}{{((n-2)!!)^{2}|S^{n-1}|}}}{Q_{n}}=\mathrm{Pfaff}+\delta _{g}B,
\label{ooo}
\end{equation}%
where $B$ is a 1-form depending locally on the metric $g$. Hence, \ by (\ref%
{gbc}), the following Gauss-Bonnet-Chern formula holds if $M$ is closed:
\begin{equation}
C_{n}\int\limits_{M}{Q_{n}\ dv_{g}}=\mathrm{Euler}(M),  \label{11.}
\end{equation}%
where $C_{n}={\frac{{1}}{{((n-2)!!)^{2}|S^{n-1}|}}}$.
\end{proposition}

Proposition~\ref{11} is one of the few results on the global properties of
the $Q_{n}$ curvature. See~\cite{B} for more details. (See also \cite{A} for
a generalization to the non-LCF case.) It establishes the integral of $Q_{n}$
as a topological quantity of a closed LCF manifolds. It is thus desirable to
extend this link to a more general class of LCF manifolds, which is one of
the motivations of this paper.

\section{A conformal Gauss-Bonnet-Chern inequality}

\label{3} \setcounter{equation}{0} \setcounter{theorem}{0}

In this Section, we focus on conformal metrics on even dimensional LCF
spaces. Assume $n=2m$ is a positive even number and $g_{0}$ is the standard
Euclidean metric on $\mathbb{R}^{n}$. A locally conformally flat metric $g$
can be represented locally as $g=g_{w}=e^{2w}g_{0}$, with $w$ being smooth.

Suppose $M$ is a manifold with a LCF metric $g$. Let $P=P_{n}$ be the
conformally invariant operator defined in (\ref{P_{2k}}). Let $Q=Q_{n}$ be
the corresponding conformal curvature invariant defined by (\ref{Q_n}). By
the discussion in Section~\ref{2}, we have
\begin{equation}
P_{w}=e^{-nw}\Delta ^{m},
\end{equation}%
\begin{equation}
Q_{g}\equiv Q_{2m,g}=e^{-nw}\Delta ^{m}w  \label{2.2}
\end{equation}%
Here without further notice the operators are all with respect to the flat
metric.

We study a LCF manifold $(M,g)$ satisfying the following assumptions:

\textbf{(A1)} $g$ is complete;

\textbf{(A2)} $R_{g}\geq 0$ near the end;

\textbf{(A3)}%
\begin{equation*}
\int_{M}{|Q}_{g}{|dv_{g}<\infty }.
\end{equation*}

Our goal of this Section is to prove Theorem~\ref{01}.

We first consider a model case, where $M=\mathbb{R}^{n}$, and the metric is
rotationally symmetric. We prove the following:

\begin{theorem}
\label{thm2} Assume $g=e^{2w}g_{0}$ is a metric on $\mathbb{R}^{n}$ such
that $w(x)=w(\Vert x\Vert )$. If $g$ satisfies assumptions {\ (A1)},{\ (A2)}
and (A3), then
\begin{equation}
C_{n}\int_{\mathbb{R}^{n}}{Q}_{g}{\ dv_{g}}\leq 1,  \label{thm2..}
\end{equation}%
with $C_{n}={\frac{1}{{((n-2)!!)^{2}|S^{n-1}|}}}$.
\end{theorem}

To prove Theorem~\ref{thm2}, we construct the Green's function as in~\cite%
{CQY1} and Theorem~\ref{thm2} is proved by solving an ODE and establishing a
uniqueness result for the conformal factor $w$.

A geometric averaging procedure is then applied to prove the following:

\begin{theorem}
\label{thm3} Assume $g=e^{2w}g_{0}$ is a metric on $\mathbb{R}^{n}$. Suppose
$g$ satisfies assumptions (A1), (A2) and (A3). Then
\begin{equation}
C_{n}\int_{\mathbb{R}^{n}}{Q}_{g}{dv_{g}}\leq 1.  \label{thm2.}
\end{equation}
\end{theorem}

Applying a gluing argument, we prove Theorem~\ref{01}, which is re-stated
below for convenience:

\begin{theorem}
\label{thm0} Assume $(M,g)$ is a LCF manifold satisfying assumptions (A1), {%
(A2)} and {(A3)}. If $M$ has only finitely many complete ends, then
\begin{equation}
C_{n}\int_{M}{Q}_{g}{dv_{g}}\leq \mathrm{Euler}(M).  \label{thm0.}
\end{equation}
\end{theorem}

This is a generalization of the Gauss-Bonnet-Chern inequality proved by
Huber~\cite{H} in the two-dimensional case and Chang, Qing and Yang~\cite%
{CQY1} in the four-dimensional case.

This Section is organized as follows. In 3.1, we prove Theorem~\ref{thm2};
in 3.2, we prove Theorem~\ref{thm3}; in 3.3, we prove Theorem~\ref{thm0}.

\subsection{$\mathbb{R}^{n}$-- the rotationally symmetric case}

Let $g=e^{2w}g_{0}$ be a conformal metric on $\mathbb{R}^{n}$ which
satisfies assumptions (A1), (A2) and (A3). In this Subsection, we make the
assumption that $w$ is rotationally symmetric; in other words, if $r=\Vert
x\Vert $, then
\begin{equation*}
w(x)=w(r).
\end{equation*}

By (\ref{2.2}), we define
\begin{equation}
f(x)\equiv Qe^{nw(x)}=\Delta ^{m}w.  \label{2.6}
\end{equation}%
Then $f\in L^{1}(\mathbb{R}^{n})$ due to (A3).

To treat the PDE~(\ref{2.6}), we notice that the standard theory of elliptic
PDE does not apply directly since the Calderon-Zygmund theory does not cover
the $L^1$ case. However, a Green's function defined in [CQY1] can still give
us a basic solution. More specifically, because $f(x)$ is integrable, the
following

\begin{equation}
v(x)\equiv C_{n}\int\limits_{\mathbb{R}^{n}}{\ln ({\frac{\Vert y\Vert }{%
\Vert x-y\Vert }})f(y)\ dy}
\end{equation}%
is well defined and smooth. It is easy to confirm that
\begin{equation}
\Delta ^{m}v(x)=f(x).
\end{equation}%
Therefore,
\begin{equation*}
\Delta ^{m}(w-v)=0.
\end{equation*}

We will now study the uniqueness for solutions for (\ref{2.6}) under the
rotational symmetry condition. (\ref{2.6}) then can be viewed as an ODE. We
prove the following simple lemma:

\begin{lemma}
\label{lemma1} If $u$ is a smooth rotationally symmetric function on $%
\mathbb{R}^{n}\backslash \{0\}$, and satisfies the differential equation
\begin{equation}
\Delta ^{m}u=0,  \label{2.1.1}
\end{equation}%
then
\begin{equation*}
u(x)=c_{0}+c_{1}\ln r+c_{2}r^{2}+c_{4}r^{4}+\cdots
+c_{n-2}r^{n-2}+c_{2}^{\prime }r^{-2}+c_{4}^{\prime }r^{-4}+\cdots
+c_{n-2}^{\prime }r^{2-n}.
\end{equation*}
\end{lemma}

\noindent \emph{Proof.} Since $u$ is rotationally symmetric, (\ref{2.1.1})
reduces to a linear ODE of n-th order. $1,\ln r,r^{2},\cdots
,r^{n-2},r^{-2},\cdots ,r^{2-n}$ are seen to be $n$ linearly independent
solutions of this ODE. Hence, the general solution is the linear combination
of these expressions. \qed

\begin{proposition}
\label{prop1} Given $f$, $v$ as above,
\begin{equation*}
\lim_{r\rightarrow 0}r\dot{v}(r)=0,
\end{equation*}%
\begin{equation*}
\lim_{r\rightarrow \infty }r\dot{v}(r)=-C_{n}\int\limits_{\mathbb{R}^{n}}{%
f(y)dy},
\end{equation*}%
where dot denotes the derivative with respect to $r$.
\end{proposition}

\noindent \emph{Proof.} Let $s\equiv \Vert y\Vert $. Clearly,
\begin{equation*}
{\frac{d}{dr}}\Vert x-y\Vert ^{2}={\frac{1}{r}}(r^{2}-s^{2}+\Vert x-y\Vert
^{2}).
\end{equation*}%
Thus,
\begin{equation*}
r\dot{v}(r)=-C_{n}\int\limits_{\mathbb{R}^{n}}{{\frac{{r^{2}-s^{2}+\Vert
x-y\Vert ^{2}}}{{2\Vert x-y\Vert ^{2}}}}f(y)dy}.
\end{equation*}

The first part of the proposition is then straightforward. The second part
is equivalent to the fact that
\begin{equation}
I(x)\equiv \int\limits_{\mathbb{R}^{n}}{{\frac{{r^{2}-s^{2}}}{{\Vert
x-y\Vert ^{2}}}}f(y)dy}\rightarrow \int\limits_{\mathbb{R}^{n}}{f(y)dy},
\label{2.14}
\end{equation}%
when $r\rightarrow \infty $.

Since $f$ is rotationally symmetric, $I$ depends only on $r=\|x\|$. Hence,

\begin{equation}
I(r)=\int\limits_{\mathbb{R}^{n}}(\int \kern-9pt\hbox{$\backslash$}_{\Vert
x\Vert =r}{{{\frac{{r^{2}-s^{2}}}{{\Vert x-y\Vert ^{2}}}}dS_{x})}f(y)dy}.
\label{2.15}
\end{equation}%
Here $dS_{x}$ is the volume form on the standard $S^{n-1}$ and
\begin{equation*}
\int \kern-9pt\hbox{$\backslash$}_{S^{n-1}}F\ dS_{x}={\frac{{%
\int_{S^{n-1}}FdS_{x}}}{{\int_{S^{n-1}}1\ dS_{x}}}}
\end{equation*}%
for a function $F$ defined on $S^{n-1}$. Define
\begin{equation}
II(r,s)\equiv \int \kern-9pt\hbox{$\backslash$}_{\Vert x\Vert =r}{\ {\frac{1%
}{\Vert x-y\Vert ^{2}}}\,dS_{x}}.  \label{2.16}
\end{equation}%
We now prove the following technical lemma:

\begin{lemma}
\label{lemma2} There exists a positive $C$ such that:
\begin{eqnarray*}
|r^{2}II(r,s)-1|\leq C|{\frac{s^{2}}{r^{2}}}| &&\mathrm{{for}\ s\leq r;} \\
II(r,s)<{\frac{C}{{s^{2}}}} &&\mathrm{{for}\ s>r.}
\end{eqnarray*}
\end{lemma}

\noindent \emph{Proof.} Taking Laplacian with respect to $y$ to $II$, we get
\begin{equation}
\Delta ^{m-2}II=C\int \kern-9pt\hbox{$\backslash$}_{\Vert x\Vert =r}{{\frac{1%
}{\Vert x-y\Vert ^{n-2}}}dS_{x}}.
\end{equation}%
Since $\frac{1}{\Vert x-y\Vert ^{n-2}}$ is a multiple of the Green's
function for the Laplacian on $\mathbb{R}^{n}$, we see that
\begin{equation}
\Delta ^{m-2}II={\frac{C}{s^{n-2}}},  \label{320}
\end{equation}%
for $s>r$;
\begin{equation}
\Delta ^{m-2}II={\frac{C}{r^{n-2}}}
\end{equation}%
for $s\leq r$, with $C$ depending only on $n$. Hence, when $s<r$, it is easy
to see from the bounded-ness of $II$ and the proof of Lemma~\ref{lemma1}
that
\begin{equation}
II(r,s)=c_{0}(r)+c_{2}(r)s^{2}+\cdots +c_{n-4}{\frac{s^{n-4}}{r^{n-2}}}.
\end{equation}

Notice the homogeneity in (\ref{2.16}), one has $r^{2}II(r,s)$ depending
only on $\frac{s}{r}$. Combining the fact that $II(r,0)={\frac{1}{r^{2}}}$,
we have
\begin{equation*}
r^{2}II(r,s)=1+p({\frac{s^{2}}{r^{2}}}),
\end{equation*}%
where $p$ is a polynomial of degree $m-1$ with no constant terms. We have
proved the lemma when $s<r$.

When $s>r$, from (\ref{320}) and the H\"{o}lder's Inequality,
\begin{equation*}
\int \kern-9pt\hbox{$\backslash$}_{\Vert x\Vert =r}{\ {\frac{1}{\Vert
x-y\Vert ^{2}}}dS_{x}}\leq (\int \kern-9pt\hbox{$\backslash$}_{\Vert x\Vert
=r}{\ {\frac{1}{\Vert x-y\Vert ^{n-2}}}dS_{x}})^{\frac{1}{{m-1}}}={\frac{C}{%
s^{2}}}.
\end{equation*}%
Lemma~\ref{lemma2} is then proved. \qed

We now continue the proof of Proposition~\ref{prop1}. Taking into account (%
\ref{2.15}), (\ref{2.16}) and Lemma~\ref{lemma2},
\begin{equation}
|I(r)-\int\limits_{\mathbb{R}^{n}}{f(y)dy}|\leq |\int\limits_{\Vert y\Vert
\leq r}{\ C{\frac{s^{2}}{r^{2}}}f(y)dy}|+C\int\limits_{\Vert y\Vert \geq r}{%
|f(y)|dy}.  \label{""}
\end{equation}

Thus, for any $\epsilon >0$, there is a positive $R$ large enough such that
\begin{equation*}
\int\limits_{s>\epsilon R}{|f(y)|dy}\leq \epsilon .
\end{equation*}%
Next, we see that when $r=\Vert x\Vert >R$,

\begin{equation*}
|\int\limits_{s\leq \epsilon r}{C{\frac{s^{2}}{r^{2}}}f(y)dy}|\leq C\epsilon
^{2}\int\limits_{\mathbb{R}^{n}}{|f(y)|dy};
\end{equation*}%
\begin{equation*}
|\int\limits_{s>\epsilon r}{C{\frac{s^{2}}{r^{2}}}f(y)dy}|\leq
C\int\limits_{s>\epsilon R}{|f(y)|dy}<C\epsilon .
\end{equation*}

Combining these and (\ref{""}), we have
\begin{equation*}
|I(x)-\int\limits_{\mathbb{R}^{n}}{f(y)dy}|\leq C(1+\int\limits_{\mathbb{R}%
^{n}}{|f(y)|dy})\epsilon ,
\end{equation*}%
the integrability of $f$ then leads to (\ref{2.14}). The proof of
Proposition~\ref{prop1} is completed. \qed

\begin{lemma}
\label{lemma3} Let $v$ as defined as above, we have, for some positive
constant $C$,
\begin{equation}
r|\dot{v}(r)|\leq C;
\end{equation}%
\begin{equation}
r^{2}|\Delta v|\leq C.
\end{equation}
\end{lemma}

\noindent \emph{Proof.} The first part follows simply from Proposition~\ref%
{prop1}. For the second part, notice that
\begin{equation}
\Delta v=C\int\limits_{\mathbb{R}^{n}}{{\frac{1}{\Vert x-y\Vert ^{2}}}f(y)dy}%
.
\end{equation}%
Again, since $f(y)$ is rotationally symmetric, we can replace $\frac{1}{%
\Vert x-y\Vert ^{2}}$ in the integrand by $II(r,s)$, which is defined in (%
\ref{2.16}). Apply Lemma~\ref{lemma2} to prove that
\begin{equation*}
r^{2}II(r,s)\leq C.
\end{equation*}%
Thus,
\begin{equation*}
|r^{2}\Delta v|\leq C\int\limits_{\mathbb{R}^{n}}{|f(y)|dy}.
\end{equation*}%
\qed

It is now possible to prove the following uniqueness result:

\begin{theorem}
\label{unique}Let the conditions be as in Theorem~\ref{thm2}. We have
\begin{equation*}
w(x)=v(x)+c,
\end{equation*}%
where $c$ is a constant.
\end{theorem}

\noindent \emph{Proof.} By Lemma~\ref{lemma1}, we have
\begin{equation*}
w(x)=v(x)+c_{0}+c_{1}\ln r+c_{2}r^{2}+c_{4}r^{4}\cdots
+c_{n-2}r^{n-2}+c_{2}^{\prime }r^{-2}+c_{4}^{\prime }r^{-4}\cdots
+c_{n-2}^{\prime }r^{2-n}.
\end{equation*}%
Since $w(x)$ is smooth at the origin, $(r\dot{w}(r))|_{r=0}=0$, by
Proposition~\ref{prop1}, we have that $c_{2k}^{\prime }=c_{1}=0$.

From (\ref{scalar-curvature}), the non-negativity of the scalar curvature is
equivalent to
\begin{equation}
\Delta w+(m-1)\Vert \nabla w\Vert ^{2}\leq 0.  \label{xx}
\end{equation}%
We prove $c_{2}=\cdots =c_{n-2}=0$ by a contradiction argument. Let $k\geq 1$
be the largest index such that $c_{2k}\not=0$. Combined with Lemma~\ref%
{lemma3}, it is clear that near infinity, $\Delta w+(m-1)\Vert \nabla w\Vert
^{2}$ has the leading term as $(m-1)c_{2k}^{2}r^{2k-2}>0$, which contradicts
with (\ref{xx}). Thus, all the $c_{i}$'s are vanishing. The proof is thus
completed. \qed

Finally, we are ready to give the following\newline
{\noindent \emph{Proof of Theorem~\ref{thm2}.} } We need to apply the
completeness condition of the metric. Notice that if
\begin{equation*}
\lim_{r\rightarrow \infty }r\dot{w}(r)=a
\end{equation*}%
exists, near $\infty $ we have that $e^{w(x)}\propto r^{a}$. For the metric $%
e^{2w(x)}$ to be complete, it has to be true that
\begin{equation}
a\geq -1.  \label{2.39}
\end{equation}%
Applying Proposition~\ref{prop1}, we have proved the inequality. \qed

From the argument above, we also have the following:

\begin{corollary}
\label{round=} Let the conditions given as in Theorem~\ref{thm2}. If further
we assume $r e^{w(x)}$ is bounded, we have the equality in (\ref{thm2..}).
\end{corollary}

\subsection{$\mathbb{R}^{n}$--the general case}

We now describe the geometric averaging procedure to reduce Theorem~\ref%
{thm3} to Theorem~\ref{thm2}.

Consider the spherical coordinate for $\mathbb{R}^{n}$. Namely, denote $x\in
\mathbb{R}^{n}$ as
\begin{equation}
x=(r,\theta ),\ \ \ r\geq 0,\ \theta \in S^{n-1},  \label{3.40}
\end{equation}%
where $S^{n-1}$ is the standard sphere (with radius $1$) in $R^{n}$.

Assume that $g=e^{2w(x)}g_{0}$ is a conformal metric on $\mathbb{R}^{n}$.
Denote $\bar{g}=e^{2\bar{w}}g_{0}$, with
\begin{equation}
\bar{w}(x)=\bar{w}(r)\equiv \int \kern-9pt\hbox{$\backslash$}_{S^{n-1}}{%
w(r,\theta )d\theta },
\end{equation}%
here we use $\int \kern-9pt\hbox{$\backslash$}_{S^{n-1}}{\cdot \,d\theta }$
to represent the average of a function over $S^{n-1}$. To study the relation
between $g$ and $\bar{g}$, denote $\nabla _{\theta }$, $\Delta _{\theta }$
as the covariant derivative and the Laplacian on $S^{n-1}$, respectively.
The following relations are obvious:
\begin{eqnarray*}
\nabla &=&\nabla _{\mathbb{R}^{n}}=(\partial _{r},{\frac{1}{r}}\nabla
_{\theta }), \\
\Delta &=&\Delta _{\mathbb{R}^{n}}=\partial _{r}^{2}+{\frac{{n-1}}{r}}%
\partial _{r}+{\frac{1}{r^{2}}}\Delta _{\theta }.
\end{eqnarray*}%
We realize that
\begin{equation*}
\Vert \nabla _{\theta }w\Vert ^{2}=r^{2}(\Vert \nabla w\Vert ^{2}-|\partial
_{r}w|^{2}).
\end{equation*}

\begin{proposition}
\label{prop3} If $g$ satisfies assumptions (A1),{\ (A2)} and{\ (A3)}, then $%
\bar{g}$ satisfies assumptions (A2) and{\ (A3)}.
\end{proposition}

\noindent \emph{Proof.} Assumption~(A3) for $\bar{g}$ metric is easy to
verify since we actually have
\begin{equation*}
\int |{Q_{n,\bar{g}}|dv_{\bar{g}}}=\int |{Q_{n,g}|dv_{g}}.
\end{equation*}%
This is because ${Q_{n,\bar{g}}dv_{\bar{g}}}=\Delta ^{m}\bar{w}\ dx$, ${%
Q_{n,g}dv_{g}}=\Delta ^{m}w\ dx$, and
\begin{equation}
\int\limits_{\Vert x\Vert =r}{\Delta ^{m}\bar{w}\ dx}=\int\limits_{\Vert
x\Vert =r}{\Delta ^{m}w\ dx}.  \label{2.50}
\end{equation}

To verify~(A2) for $\bar{g}$, by (\ref{scalar-curvature}), $R_{g}\geq 0$ is
equivalent to $\Delta w+(m-1)\Vert \nabla w\Vert ^{2}\leq 0$. Since
\begin{equation*}
\Delta \bar{w}=\int \kern-9pt\hbox{$\backslash$}{\Delta wd\theta },
\end{equation*}%
and
\begin{equation*}
\Vert \nabla \bar{w}\Vert ^{2}=(\int \kern-9pt\hbox{$\backslash$}\partial
_{r}wd\theta )^{2}\leq \int \kern-9pt\hbox{$\backslash$}{\Vert \nabla
w(r,\theta )\Vert ^{2}d\theta }.
\end{equation*}%
Hence it is apparent that $\Delta \bar{w}+(m-1)\Vert \nabla \bar{w}\Vert
^{2}\leq 0$, which implies that $R_{\bar{g}}$ is non-negative.\qed {}

\begin{proposition}
\label{complete}If $g$ satisfies assumptions (A1),{\ (A2)} and{\ (A3)}, then
$\bar{g}$ is complete.
\end{proposition}

\noindent \emph{Proof.} As before, we define \
\begin{eqnarray}
f(x) &\equiv &Qe^{nw(x)}=\Delta ^{m}w,  \notag \\
v(x) &\equiv &C_{n}\int\limits_{\mathbb{R}^{n}}{\ln ({\frac{\Vert y\Vert }{%
\Vert x-y\Vert }})f(y)\ dy} \\
u(x) &=&w(x)-v(x).
\end{eqnarray}%
Then $f\in L^{1}(\mathbb{R}^{n})$ and $v,u\in C^{\infty }(\mathbb{R}^{n})$.

We first show two intermediate results.\ The first one is a generalization
of Theorem \ref{unique}:

\begin{claim}
\label{claim1}$u(x)$ is a constant function.
\end{claim}

\begin{proof}
Similar to the rotationally symmetric case, we have
\begin{equation*}
\Delta ^{m}u=0.
\end{equation*}%
To show the uniqueness result, we proceed to consider the rotational
symmetrization procedure with respect to a fixed point $P\in \mathbb{R}^{n}:$%
\begin{eqnarray*}
\bar{w}_{P}(x) &=&\dint \kern-9pt\hbox{$\backslash$}_{||y-P||=\left\vert
\left\vert x-P\right\vert \right\vert }w(y)d\theta \\
\bar{v}_{P}(x) &=&\dint \kern-9pt\hbox{$\backslash$}_{||y-P||=\left\vert
\left\vert x-P\right\vert \right\vert }v(y)d\theta \\
\bar{u}_{P}(x) &=&\dint \kern-9pt\hbox{$\backslash$}_{||y-P||=\left\vert
\left\vert x-P\right\vert \right\vert }u(y)d\theta .
\end{eqnarray*}%
Notice that
\begin{eqnarray}
v(x) &=&C_{n}\int\limits_{\mathbb{R}^{n}}{\ln ({\frac{\Vert y-P\Vert }{\Vert
x-y\Vert }})f(y)\ dy+}C_{n}\int\limits_{\mathbb{R}^{n}}{\ln ({\frac{\Vert
y\Vert }{\Vert P-y\Vert }})f(y)\ dy}  \label{ah} \\
&=&C_{n}\int\limits_{\mathbb{R}^{n}}{\ln ({\frac{\Vert y-P\Vert }{\Vert
x-y\Vert }})f(y)\ dy+C,}  \notag
\end{eqnarray}%
apply Proposition \ref{prop3}, $g_{\bar{w}_{P}}$ satisfies conditions (A2)
and (A3). Apply Lemma~\ref{lemma1} and Theorem \ref{unique}, we have
\begin{equation*}
\bar{u}_{P}(x)=c.
\end{equation*}%
In particular, it implies that
\begin{equation*}
\Delta u(P)=\Delta \bar{u}_{P}(P)=0.
\end{equation*}%
Hence, we have shown that $u$ is harmonic over $\mathbb{R}^{n}.$ We can
finish the prove of Claim~\ref{claim1} by following an argument of~\cite%
{CQY1}: because $u$ is harmonic, so is $u_{i}(x)=\frac{\partial u}{\partial
x^{i}}(x).$ It leads to
\begin{eqnarray}
|u_{i}(P)|^{2} &=&\left\vert \int \kern-9pt\hbox{$\backslash$}%
_{||x-P||=r}u_{i}d\theta \right\vert ^{2}  \notag \\
&=&\left\vert \int \kern-9pt\hbox{$\backslash$}_{||x-P||=r}u_{i}d\theta
\right\vert ^{2}  \notag \\
&\leq &\int \kern-9pt\hbox{$\backslash$}_{||x-P||=r}\left\vert \left\vert
\nabla u\right\vert \right\vert ^{2}d\theta  \notag \\
&\leq &\int \kern-9pt\hbox{$\backslash$}_{||x-P||=r}(\left\vert \left\vert
\nabla w\right\vert \right\vert ^{2}+\left\vert \left\vert \nabla
v\right\vert \right\vert ^{2})d\theta \leq \frac{C}{r^{2}}\rightarrow 0,
\end{eqnarray}%
as $r\rightarrow \infty .$ In the last step we have applied Lemma \ref%
{lemma3} for (\ref{ah}) and the fact that

\begin{eqnarray*}
&&\int \kern-9pt\hbox{$\backslash$}_{||x-P||=r}\left\vert \left\vert \nabla
w\right\vert \right\vert ^{2}d\theta \\
&=&\frac{1}{m-1}\int \kern-9pt\hbox{$\backslash$}_{||x-P||=r}(\Delta
w-e^{w}J_{w})d\theta \\
&\leq &\frac{1}{m-1}\Delta \bar{w}_{P}(r)=0.
\end{eqnarray*}%
We thus have proved that $u_{i}=0$ all any $x^{i};$ hence, $u$ is a constant.
\end{proof}

The second intermediate result is the following analogue of Lemma 3.2 of
\cite{CQY1}:

\begin{claim}
\label{claim2}If
\begin{equation*}
w(x)=C_{n}\int\limits_{\mathbb{R}^{n}}{\ln ({\frac{\Vert y\Vert }{\Vert
x-y\Vert }})\Delta }^{m}w{(y)\ dy+C,}
\end{equation*}%
then
\begin{equation*}
e^{-\bar{w}}\dint \kern-9pt\hbox{$\backslash$}_{||x||=r}e^{w}\rightarrow 1,
\end{equation*}%
as $r\rightarrow \infty .$
\end{claim}

The proof of Claim~\ref{claim2} is identical to the proof of Lemma 3.2 of
\cite{CQY1}, which treats dimension four case. We omit it here.

Now we can continue the proof of Proposition \ref{complete}. We only need to
show that $\int_{0}^{\infty }e^{\bar{w}}dr$ is divergent.\ Since for a fixed
$\theta ,$ $\int_{0}^{\infty }e^{w(r,\theta )}dr$ is divergent because of
completeness of metric $g,$ this can be proved by applying Claims \ref%
{claim1} and \ref{claim2}.

Thus we have completed the proof of Proposition \ref{complete}.\qed

\bigskip

\begin{corollary}
\label{coro1} If $g=e^{2w}g_{0}$ is a conformal metric on $\mathbb{R}^{n}$
such that $R_{g}\geq 0$, and $P_{n}=0$, then $w$ is a constant.
\end{corollary}

\noindent \emph{Proof.} One constructs the metric $\bar g=e^{2\bar w}$ as in
(\ref{3.40}). From the proof of Proposition~\ref{prop3}, it is true that $%
\bar g$ has non-negative scalar curvature and vanishing $Q$ curvature. The
conclusion then follows from Theorem~\ref{unique}. \qed

We now give the proof of Theorem~\ref{thm3}. But this is a straightforward
application of Propositions~\ref{prop3}, \ref{complete} and Theorem~\ref%
{thm2}. \qed

\subsection{LCF manifolds with finitely many ends}

In this Subsection, we give the proof of Theorem~\ref{thm0}, which is an
extension of Theorem~1.2 of [CQY2] in higher dimensional case. We will take
advantage of the topological invariance of $\int Qdv$ and give a doubling
argument. Our approach is more geometrical, comparing to the approach of\
Chang-Qing-Yang, which is more analytical.

First we prove the following simplified result:

\begin{proposition}
\label{hmm1} Assume $\Omega$ is a domain in $S^n$ with a conformal metric $g$
satisfying assumptions~(A1), (A2) and (A3). If $\Lambda=S^n\backslash\Omega$
is a finite set of $k$ points, then
\begin{equation}
C_n\int_\Omega{Q\ dv_g}\leq (2-k).  \label{hmm1.}
\end{equation}
\end{proposition}

\noindent \emph{Proof.} Let $\Lambda =\{p_{1},\cdots ,p_{k}\}$. A
stereographic projection from $S^{n}$ to $\mathbb{R}^{n}$ can be chosen so
that $p_{1}$ is sent to infinity. Without confusion, we identify the images
of $\Lambda $ under the projection with itself. There is a function $w$
smooth away from $\Lambda $ such that the metric can be represented as $%
g=e^{2w}g_{0}$, where $g_{0}$ is the Euclidean metric on $\mathbb{R}^{n}$.
We fix a partition of unity,
\begin{equation*}
1=l_{1}(x)+\cdots +\l _{k}(x),
\end{equation*}%
such that $l_{i}(x)$ is a smooth function supported near $p_{i}$ and $%
l_{i}=1 $ near $p_{i}$. Let $w_{i}(x)=w(x)l_{i}(x)$. We consider the metric $%
g_{i}=e^{2w_{i}}g_{0}$.

Note that $g_{1}$ satisfies assumptions (A1), (A2) and {(A3)}. Theorem~\ref%
{thm3} then gives the follows:
\begin{equation}
C_{n}\int\limits_{\mathbb{R}^{n}}{Q_{g_{1}}\ dv_{g_{1}}}=C_{n}\int\limits_{%
\mathbb{R}^{n}}{\Delta ^{m}(w_{1})dx}\leq 1.  \label{point1}
\end{equation}

For a fixed $i\geq 2$, without loss of generality, we assume that $p_{i}$ is
just the origin. $w_{i}$ has compact support and the metric $g_{k}$ also
satisfies assumptions (A1), (A2) and {(A3)}. We construct
\begin{equation*}
\bar{w}_{i}\equiv \int \kern-9pt\hbox{$\backslash$}_{S^{n-1}}{w_{i}(r,\theta
)\ d\theta }
\end{equation*}%
and
\begin{equation*}
v_{i}(x)\equiv C_{n}\int\limits_{\mathbb{R}^{n}}{\ln ({\frac{\Vert y\Vert }{%
\Vert x-y\Vert }})\Delta ^{m}\bar{w}_{i}(y)\ dy}.
\end{equation*}%
Notice that Proposition~\ref{prop1} and Proposition~\ref{prop3} can still be
applied to the metric $g_{i}$. Tracing the argument in the proof of Theorem~%
\ref{unique}, we see that
\begin{equation*}
\bar{w}_{i}(x)=v_{i}(x)+c_{1,i}\ln r+c_{0,i}.
\end{equation*}%
Note that $\bar{w}=0$ for $\Vert x\Vert $ large, by Proposition~\ref{prop1},
\begin{equation}
c_{1,i}=C_{n}\int\limits_{\mathbb{R}^{n}}{Q(\bar{g}_{i})dv_{\bar{g}_{i}}}%
=C_{n}\int\limits_{\mathbb{R}^{n}}\Delta ^{m}\bar{w}\ dx=C_{n}\int\limits_{%
\mathbb{R}^{n}}{\Delta ^{m}w_{i}\ dx}.  \label{qqq}
\end{equation}

Follow the proof of Theorem~\ref{thm2}, instead of getting (\ref{2.39}), the
completeness of $\bar{g}_i$ near the origin shows
\begin{equation}
c_{1,i}\leq -1.  \label{pointk}
\end{equation}

Combine (\ref{2.50}), (\ref{point1}), (\ref{qqq}), (\ref{pointk}) with the
fact that $w=\sum w_{k}$, we have
\begin{equation*}
C_{n}\int\limits_{\Omega }{Q_{g}\ dv_{g}}=C_{n}\int\limits_{\mathbb{R}^{n}}{%
\Delta ^{m}w\,dx}=\sum_{i}C_{n}\int\limits_{\mathbb{R}^{n}}{\Delta
^{m}w_{i}\,dx}\leq 1+(k-1)(-1)=2-k.
\end{equation*}%
Thus, Proposition~\ref{hmm1} is proven. \qed

Define $d(\cdot)$ to be the distance function to $\Lambda$ on $S^{n-1}$. We
have the following easy extension of Corollary~\ref{round=}.

\begin{corollary}
\label{finite=} Let the conditions be those of Proposition~\ref{hmm1}. If
the conformal factor $w$ satisfies that $d(p)e^w(p)\leq C$ for any $p\in
\Omega$ and some positive constant $C$, then the equality holds in (\ref%
{hmm1.}).
\end{corollary}

We are ready to give the

{\noindent \emph{Proof of Theorem~\ref{thm0}.} } Suppose a complete LCF
manifold $M$ has $k$ disjoint ends $E_{1},\cdots ,E_{k}$. We can thus choose
a local coordinate chart for each $E_{i}$ such that the metric is
represented as $e^{2w_{i}(x)}g_{0}$, $g_{0}$ being the $n$-dimensional
Euclidean metric and $\Vert x\Vert >1$. Our approach is following: First we
modify the conformal metric so that each end links the the manifold $M$ in a
strict tubular fashion; then, we cut off each ends to get manifolds with
boudary; finally, in order to esitmate $\int Qdv,$ we double the compact
piece and apply Gauss-Bonnet-Chern formula and extend each ends naturally to
apply Proposition~\ref{hmm1}.

First, we do a compact perturbation of the metric near the ends. Let $\eta $
be a cut off function such that $\eta (x)=1$ for $2<\Vert x\Vert <3$ and $%
\eta (x)=0$ for $\Vert x\Vert <1$ and $\Vert x\Vert >4$. Define a new
conformal metric
\begin{equation*}
g^{\prime }=\left\{
\begin{array}{l}
e^{2\eta (x)(-w_{i}(x)-\ln (\Vert x\Vert ))}g\ \ \ \ \mathrm{{on}\ E_{i};}
\\
g\ \ \ \ \ \ \ \ \ \ \ \ \ \ \ \ \ \ \ \ \ \ \ \ \ \ \ \ \mathrm{{on}\
M-\cup _{i}E_{i}.}%
\end{array}%
\right.
\end{equation*}%
By the choice of $\eta $, $g^{\prime }$ is well-defined and smooth. Hence,
\begin{eqnarray}
&&\int_{M}{Q_{g^{\prime }}\ dv_{g^{\prime }}}-\int_{M}{Q_{g}\ dv_{g}}
\label{g'} \\
&=&\sum_{i}\int_{\o {1}{2}\leq \Vert x\Vert \leq {\o {9}{2}}}\Delta
^{m}[\eta (x)(-w_{i}(x)-\ln \Vert x\Vert ]\ dv_{g_{0}}=0.  \notag
\end{eqnarray}

Second, for each $i$, define $E_{i}^{\prime }\subset E_{i}=\{x;\Vert x\Vert >%
{\frac{5}{2}}\}$. If $M_{1}\equiv M-\cup _{i}E_{i}^{\prime }$, then $M_{1}$
is a compact manifold with boundary. $\partial M_{1}$ has $k$ components and
near each of them the metric $g^{\prime }$ is a locally product metric due
to the construction of $g^{\prime }$. Hence, we can glue two pieces of $%
(M_{1},g^{\prime })$ together to get a closed manifold $M_{2}$. Still
referring to the gluing metric on $M_{2}$ as $g^{\prime }$, we apply the
Gauss-Bonnet-Chern formula for closed manifolds to get
\begin{equation}
C_{n}\int_{M_{1}}{Q(g^{\prime })\ dv_{g^{\prime }}}=\o {1}{2}%
C_{n}\int_{M_{2}}{Q(g^{\prime })\ dv_{g^{\prime }}}=\o {1}{2}\mathrm{Euler}%
(M_{2})=\mathrm{Euler}(M_{1}).  \label{M_1}
\end{equation}

Next, the metric $g^{\prime}$ on $E^{\prime}_i$ is equal to ${\frac{1}{{%
\|x\|^2}}}g_0$ near the boundary of $E^{\prime}_i$, can be extended in the
coordinate chart to region $E^{\prime\prime}_i=\{ \|x\|>0\}$, still denoted
as $g^{\prime}$, so that $g^{\prime}={\frac{1}{{\|x\|^2}}}g_0$ for $x\in
E^{\prime\prime}_i-E^{\prime}_i$. Notice for $x\in
E^{\prime\prime}-E^{\prime}$, $Q(g^{\prime})=\|x\|^n \Delta^m {(-\ln \|x\|)}%
=0$. $(E^{\prime\prime}_i, g^{\prime})$ satisfies assumptions {(A1)}, (A2)
and {(A3)}. Proposition~\ref{hmm1} is applied to $(E^{\prime\prime}_i,
g^{\prime})$ to get
\begin{equation}
\int_{E^{\prime}_i}
Q(g^{\prime})dv_{g^{\prime}}=\int_{E^{\prime\prime}_i}Q(g^{\prime})dv_{g^{%
\prime}}\leq 0 .  \label{E'}
\end{equation}

Finally, combine (\ref{g'}), (\ref{M_1}) and (\ref{E'}), we prove Theorem~%
\ref{thm0}.\qed

\section{Conformal compactification for certain LCF 4-folds with integrable
Pfaffian curvature}

\label{5} \setcounter{equation}{0} \setcounter{theorem}{0}

In the previous Sections, we proved the Gauss-Bonnet-Chern-type inequality
for certain LCF manifolds with integrable $Q$ curvature. As we have seen,
the conformal variation of $Q$ is just $P$ operator, which is linear
elliptic. This fact is crucial in our study in Section~\ref{2}. However, it
is interesting that some of the techniques developed in the previous
Sections are also applicable to to study the Pfaffian curvature of certain
complete LCF manifolds of dimension four. As in~[CQY2], we pose the
following stronger assumption for the curvature:
\begin{equation}
C_{1}\geq R_{g}\geq C>0,\ \Vert \nabla _{g}R_{g}\Vert _{g}\leq C,\ \mathrm{%
Ric}_{g}\geq -Cg  \tag{A4}  \label{A4}
\end{equation}%
for some positive constants $C$ and $C_{1}$.

\bigskip

In this Section, we prove Theorem \ref{03}, which is rephrased here for
readers' convenience:

\begin{theorem}
\label{4-1} Let $\Omega $ be a manifold with a LCF metric satisfying
assumptions (A1), (A2), and (A4). If further
\begin{equation*}
\int\limits_{\Omega }{|\mathrm{Pfaff}_{g}|dv_{g}}<\infty ,
\end{equation*}%
then $\Lambda =S^{4}\backslash \Omega $ is a finite set.
\end{theorem}

Comparing to the situation treated in [CQY2], we replace $Q$ curvature by
the Pfaffian of the manifold. Hence, it is a non-linear extension of the
main result of [CQY2].

By the curvature condition and the extension map construction of Schoen-Yau~%
\cite{SY}, we can view $\Omega $ as a domain in $S^{4}.$As in before, we
identify $\Omega $ with its image in $\mathbb{R}^{n}$ under a stereographic
projection and write $g=e^{2w}g_{0}$ with $g_{0}$ being the Euclidean
metric. Again, we use the upper-bar to denote geometric quantities with
respect to $g$ metric.

First we study the fully non-linear transformation law of the Pfaffian. For $%
g$, by (\ref{pan4}) and (\ref{scalar-curvature}), and (\ref{ricci}) and (\ref%
{ooo}), we have the following:
\begin{equation}
\mathrm{Pfaff}_{g}=C_{4}(Q_{g}+\bar{\Delta}J).  \label{faf1}
\end{equation}%
\begin{equation}
\mathrm{Pfaff}=C_{4}e^{-4w}((\Delta w)^{2}-\Vert \nabla ^{2}w\Vert
^{2}+2\nabla ^{2}w(\nabla w,\nabla w)+\Vert \nabla w\Vert ^{2}\Delta w),
\label{faf3}
\end{equation}%
where all the operators are with respect to the Euclidean metric and $\nabla
^{2}$ denotes the Hessian.

It is an interesting observation that
\begin{equation}
\mathrm{Pfaff}=C_{4}\sigma _{2},  \label{faf2}
\end{equation}%
where $\sigma _{2}$ is defined in (\ref{sigma}). This is actually a special
case of a more general fact that for any LCF metric on a $2m$-dimensional
manifold, the Pfaffian is a constant multiple of $\sigma _{m}$ (Cf.~\cite{V}%
, for example).

In this section, we first give a $C^{0}$ estimate of the conformal factor $%
w; $ then we give an estimate of the size of $S^{4}\backslash \Omega .$

Notice that in this Section, we fix $n=4$ and $m=2,$ though many arguments
work for general dimensions. See~\cite{F} for more general statements.

\subsection{$\mathbf{C}^{0}$ estimate}

\label{4001}

In this Subsection, we give the key estimates of the conformal factor.

First we quote a lemma of Yau (Cf.~\cite{SY}), which is a special case of
the gradient estimate for positive harmonic functions on a complete manifold.

\begin{lemma}
\label{savior} For a manifold $M$ with a LCF complete metric $g=e^{2w}g_{0}$
satisfying the following: the scalar curvature $R_{g}$ and the Ricci
curvature $\mathrm{Ric}_{g}$ satisfy the following point-wise estimates near
the complete end:
\begin{equation*}
C\geq R_{g}\geq 0,\ \Vert \nabla _{g}R_{g}\Vert _{g}\leq C,\ \mathrm{Ric}%
_{g}\geq -Cg
\end{equation*}%
for some positive constant $C$, then there exists a constant $C$ such that
\begin{equation*}
\Vert \nabla _{g}w\Vert _{g}\leq C.
\end{equation*}
\end{lemma}

We then prove a non-existence result which is an analogue of Theorem~\ref%
{unique}.

\begin{lemma}
\label{pf-unique} There is no metric $g=e^{2w}g_{0}$ on $\mathbb{R}^{4}$
satisfying (A1), (A2), (A4) and
\begin{equation}
\mathrm{Pfaff}=0.  \label{three}
\end{equation}
\end{lemma}

Notice now that the averaging method we applied to prove Theorem~\ref{unique}
does not work since unlike the Q curvature, the Pfaffian does not satisfy a
linear transformation law. An integral estimate is applied instead.

\noindent \emph{Proof.} We assume there exists such a metric. First, we
prove that the Ricci curvature is bounded. Because $0=\sigma _{2}=\o {1}{2}%
(-\Vert A\Vert ^{2}+J^{2})$, by (\ref{A_g}) and (A2),
\begin{equation}
\Vert \mathrm{Ric}_{g}\Vert ^{2}=\Vert (n-2)A_{g}+J\Vert ^{2}\leq
C(R_{g}^{2}+\Vert A_{g}\Vert ^{2})\leq C^{\prime }  \label{four}
\end{equation}%
for some positive constants $C$ and $C^{\prime }$.

It is clear that Lemma~\ref{savior} is applicable. Thus, for some positive $%
C $,
\begin{equation}
\|\nabla w\|\leq Ce^w.  \label{five}
\end{equation}

For any region $D\subset \Omega $ and a positive $\alpha <\o {1}{2}$, we
claim that there is some positive constant $C$ such that
\begin{equation}
\alpha \int\limits_{D}{e^{(4+\alpha )w}dv}\leq C\int\limits_{\partial D}{%
e^{(3+\alpha )w}dv^{\prime }},  \label{zero}
\end{equation}%
with $dv$, $dv^{\prime }$ denoting the Euclidean volume forms on $D$ and $%
\partial D$, respectively.

We now prove the claim. Notice that by (\ref{faf3}) and (\ref{faf2}), $%
e^{4w}\sigma _{2}$ can be re-written as a divergence form:
\begin{equation}
e^{4w}\sigma _{2}=\delta ((\Delta w+\Vert \nabla w\Vert ^{2}-\nabla
^{2}w)dw).
\end{equation}%
It follows that
\begin{eqnarray}
\int\limits_{D}{\sigma _{2}e^{(4+\alpha )w}dv} &=&\int\limits_{D}{(-\Delta
w-\Vert \nabla w\Vert ^{2}+\nabla ^{2}w)(dw,\alpha e^{\alpha w}dw)dv} \\
&&+\int\limits_{\partial D}{e^{\alpha w}(\Delta w+\Vert \nabla w\Vert
^{2}-\nabla ^{2}w)(\partial _{\mathbf{n}}w)dv^{\prime }}.  \notag
\end{eqnarray}

Since $\sigma _{2}=\mathrm{Pfaff}=0$,
\begin{eqnarray}
0 &=&\int\limits_{D}{\alpha e^{\alpha w}(-\Delta w-\Vert \nabla w\Vert
^{2})\Vert \nabla w\Vert ^{2}-\o {\Vert \nabla w\Vert ^{2}}{2}(\alpha
e^{\alpha w}w_{i})_{i}dv}  \label{0.5} \\
&&+\int\limits_{\partial D}{e^{\alpha w}[(\Delta w+\Vert \nabla w\Vert
^{2}-\nabla ^{2}w)(\partial _{\mathbf{n}}w)+}\frac{{{\Vert \nabla w\Vert ^{2}%
}}}{2}{(\alpha \partial _{\mathbf{n}}w)]}  \notag \\
&=&\int\limits_{D}{\alpha e^{\alpha w}\Vert \nabla w\Vert ^{2}[-\Delta
w-\Vert \nabla w\Vert ^{2}-}\frac{{\alpha }}{2}{\Vert \nabla w\Vert
^{2}-\Delta w]dv}  \notag \\
&&\ \ +\int\limits_{\partial D}{e^{\alpha w}[(\Delta w+\Vert \nabla w\Vert
^{2}-\nabla ^{2}w)(\partial _{\mathbf{n}}w)+\o {\Vert \nabla w\Vert ^{2}}{2}%
(\alpha \partial _{\mathbf{n}}w)]dv^{\prime }}.  \notag
\end{eqnarray}

By (A2),
\begin{equation}
|\Delta w+ \|\nabla w\|^2|= |J|e^{2w}\leq Ce^{2w},  \label{six}
\end{equation}
for some constant $C$. Combine (\ref{six}) and (\ref{five}) we get
\begin{equation}
|\Delta w|\leq C e^{2w}  \label{seven}
\end{equation}
Combine (\ref{four}), (\ref{five}) and (\ref{seven}), it is not hard to see
that
\begin{equation}
\|\nabla^2 w\|\leq C e^{2w}  \label{eight}
\end{equation}
Therefore, from (\ref{five}), (\ref{seven}) and (\ref{eight}), we have
\begin{equation}
\int\limits_{\partial D}{e^{\alpha w} [(\Delta w+\|\nabla w\|^2-\nabla^2
w)(\partial_{\mathbf{n}} w)+\o {\|\nabla w\|^2}{2} (\alpha \partial_{\mathbf{%
n}} w)]} \leq C\int\limits_{\partial D}{e^{(3+\alpha)w}dv^{\prime}}.
\label{8.5}
\end{equation}

On the other hand,
\begin{equation*}
\int\limits_{D}{\ e^{\alpha w}\Vert \nabla w\Vert ^{2}[-\Delta w-\Vert
\nabla w\Vert ^{2}-\o {\alpha }2\Vert \nabla w\Vert ^{2}-\Delta w]dv}\geq
C\int\limits_{D}{e^{\alpha w}\Vert \nabla w\Vert ^{2}Je^{2w}dv}
\end{equation*}%
\begin{equation}
\geq C\int\limits_{D}{e^{\alpha w}\Vert \nabla (e^{w})\Vert ^{2}dv}.
\label{nine}
\end{equation}%
Since
\begin{equation*}
J=-e^{-3w}\Delta e^{w}
\end{equation*}
in dimension four, through integration by part,
\begin{equation*}
\int\limits_{D}{e^{\alpha w}\Vert \nabla (e^{w})\Vert ^{2}dv}=\int\limits_{D}%
{(e^{\alpha w}Je^{4w}-\alpha e^{\alpha w}\Vert \nabla (e^{w})\Vert ^{2})dv}%
+\int\limits_{\partial D}{e^{(2+\alpha )w}\partial _{\mathbf{n}}w\
dv^{\prime }}.
\end{equation*}%
With $\alpha \leq \frac{\text{1}}{2}$, it is true that
\begin{equation}
\int\limits_{D}{e^{\alpha w}\Vert \nabla (e^{w})\Vert ^{2}dv}\leq
C\int\limits_{D}{e^{(4+\alpha )w}dv}-\alpha C\int\limits_{\partial D}{%
e^{(3+\alpha )w}dv^{\prime }}.  \label{ten}
\end{equation}%
Combine (\ref{0.5}), (\ref{ten}),(\ref{nine}) and (\ref{8.5}), we reach the
proof of the claim (\ref{zero}).

Choose the domain $D$ as $B(0,r)$, the ball centered at the origin with
radius $r$, and define
\begin{equation}
F(r)\equiv\int\limits_{B(0,r)}{e^{(4+\alpha)w}dv}.  \label{10.1}
\end{equation}
It is clear that
\begin{equation}
F^{\prime}(r)=\int\limits_{\partial B(0,r)} {e^{(4+\alpha)w}dv^{\prime}}.
\label{10.2}
\end{equation}

By H\"{o}lder Inequality,
\begin{equation}
\int\limits_{\partial D}{e^{(3+\alpha)w}dv^{\prime}}\leq
[\int\limits_{\partial D}{e^{(4+\alpha)w}dv^{\prime}}]^{\o {3+\alpha}{%
4+\alpha}}[\int\limits_{\partial D}{dv^{\prime}}]^{\o {1}{4+\alpha}}.
\label{eleven}
\end{equation}
Substitute (\ref{10.1}), (\ref{10.2} and (\ref{eleven}) into (\ref{zero}),
we have that
\begin{equation}
\alpha F(r)\leq C [F^{\prime}(r)^{\o {3+\alpha}{4+\alpha}}]\cdot r^{\o {3}{%
4+\alpha}},
\end{equation}
which implies that
\begin{equation}
(-[F(r)]^{\o {-1}{3+\alpha}})^{\prime}\geq C\alpha^{\o {7+2\alpha}{3+\alpha}%
}(r^{\o {\alpha}{3+\alpha}})^{\prime}.  \label{twelve}
\end{equation}

Then for any $b>1$, integrate (\ref{twelve}) over $[1,b]$ we get
\begin{equation*}
\lbrack F(b)]^{\o {-1}{3+\alpha }}\leq C\alpha ^{\o {7+2\alpha }{3+\alpha }%
}(1-b^{\o {\alpha }{3+\alpha }})+[F(1)]^{\o {-1}{3+\alpha }}.
\end{equation*}%
Let $b$ tends to $\infty $ we would get the absurd conclusion that $F(b)\leq
0$. We thus have finished the proof of Lemma~\ref{pf-unique}.

Now we are ready to show the $C^1$ estimate as follows:

\begin{lemma}
\label{blow-up} For a LCF metric $g$ on $\Omega \subset S^{n}$ satisfying
(A1), (A2), and (A4), we have
\begin{equation}
C<e^{w(x)}d(x)\leq C^{\prime },
\end{equation}%
for some positive $C$ and $C^{\prime }$.
\end{lemma}

\noindent \emph{Proof.} The left-hand side of the inequality is a direct
consequence of Lemma~\ref{savior}.

To prove the right-hand side of the inequality, we run a blow-up argument,
following Schoen [S] and [CQY2]. For simplicity, denote
\begin{equation}
u(x)=e^{(m-1)w(x)}.
\end{equation}%
If the claim is not true, we would have a sequence of $\{x_{i}\}\in \Omega $%
, such that
\begin{equation}
A_{i}=u(x_{i})d^{m/2-1}(x_{i})\rightarrow \infty .  \label{3.3}
\end{equation}

For simplicity, we define the following quantities:
\begin{equation}
\sigma_i\equiv {\frac{1}{2}} d(x_i),\ \ f_i(y)\equiv (\sigma_i-d(y,x_i))^{{%
m/2}-1} u(y).
\end{equation}
It follows from (\ref{3.3}) that
\begin{equation}
f_i(x_i)=\sigma_i^{m/2-1} u(x_i)={\frac{1}{2}} A_i\to \infty
\end{equation}
when $i \to \infty$ and $f_i (y)=0$ for $y\in \partial B(x_i,\sigma_i)$.
Thus, there exists some point $y_i$ such that
\begin{equation}
f(y_i)=\max \{f_i(y):y\in B(x_i,\sigma_i)\}.
\end{equation}
Set
\begin{equation}
\lambda_i\equiv u(y_i),
\end{equation}
and
\begin{equation}
v_i(x)\equiv \lambda_i^{-1} u(\lambda_i^{-{\frac{1}{{m-1}}} } x+y_i)
\end{equation}
for $x\in B(0,R_i)$. Hence,
\begin{equation}
v_i(0)=1.
\end{equation}

Let $r_i=\o {1}{2} (\sigma_i-d(x_i, y_i)$ and $R_i=r_i u(y_i)$; then $x\in
B(0, R_i)$ if and only $y=\lambda_i^{-{\frac{1}{{m-1}}}} x+y_i\in B(x_i,
r_i) $. Notice that $R_i\to \infty$ as $i\to \infty$.

Notice that
\begin{equation}
0<v_i(x)={\frac{{u(y)}}{{u(y_i)}}}\leq ({\frac{{\sigma_i-d(x_i, y_i)}}{{%
\sigma_i-d(x_i,y)}}})^{m/2-1}\leq ({\frac{\sigma_i}{{\sigma_i-r_i}}}
)^{m/2-1}\leq 2^{m/2-1}.
\end{equation}

By (\ref{scalar-curvature}), $g_{v_i}$ satisfies the following
\begin{equation}
J_i(x)=-\Delta v_i(x) v_i(x)^{-{\frac{{n+2}}{{n-2}}}},
\end{equation}
for $x\in B(0,R_i)$. Note that the bounded-ness of $v_i$ will also give the
bounded-ness of $|\nabla_x J_i(x)|$.

It thus follows that, taking a subsequence if necessary,
\begin{equation}
J_{i}\rightarrow J_{\infty }\in C_{loc}^{\alpha }(\mathbb{R^{n}})
\end{equation}%
for some $J_{\infty }\geq C>0$. Hence a subsequence of $v_{i}$ converges
uniformly on compact sets in $C^{1,\alpha }(\mathbb{R}^{n})$. Let the limit
function be $v_{\infty }$. By the standard elliptic theory, we show that $%
v_{\infty }\in C^{2,\alpha }(\mathbb{R}^{n})$, and
\begin{equation}
-\Delta v_{\infty }(x)=J_{\infty }(x)v_{\infty }(x)^{{\frac{{m+1}}{{m-1}}}}.
\end{equation}%
Applying the elliptic theory again, $v_{\infty }$ is actually smooth.
Applying the maximum principle and the fact that $v_{\infty }(0)=1$, we
derive $v_{\infty }(x)>0$. If $w_{i}(x)\equiv \ln v_{i}(x)$ and $w_{\infty
}(x)\equiv \ln v_{\infty }(x)$, by passing to a subsequence, we conclude
\begin{equation}
w_{i}\rightarrow w_{\infty }\ \ \mathrm{{in}\ C^{2,\alpha }(\mathbb{R}}^{n}%
\mathrm{).}
\end{equation}%
This implies that for $g_{\infty }=e^{2w_\infty }g_{0},$%
\begin{equation*}
\text{Pfaff}_{g_{\infty }}=0.
\end{equation*}%
It is easy to see that $g_{\infty }$ satisfies the condition of Lemma \ref%
{pf-unique}, hence we have a contradiction.\ We have thus completed the
proof of Lemma~\ref{blow-up}.

\qed

\begin{remark}
By the way we present the proof of Lemma~\ref{blow-up}, the blowup argument
we used works also for general dimensions.
\end{remark}

\subsection{Proof of Theorem~\protect\ref{4-1}}

Following [CQY2], we define subsets of $M:$
\begin{eqnarray*}
U_{\lambda } &=&\{x:\ e^{w(x)}\geq \lambda \}, \\
S_{\lambda } &=&\{x:\ e^{w(x)}=\lambda \}.
\end{eqnarray*}%
Then, if $\mathbf{n}$ is the outward normal vector of $S_{\lambda }$ as the
boundary of $U_{\lambda }$, $\partial _{\mathbf{n}}w\geq 0$. We work with
the level sets of $w$ from now on.

We begin with a technical result:

\begin{lemma}
\label{reduce} Given any LCF metric $g$ on $\Omega $, and any $f\in
C^{\infty }(\Omega )$, if $f=f(\lambda )$, then
\begin{equation}
\int\limits_{U_{\lambda }}{(\Delta _{g}f)dv_{g}}=\lambda {\frac{d}{{d\lambda
}}}[\int\limits_{U_{\lambda }}{(\Delta _{g}w)f\ dv_{g}}-\int\limits_{S_{%
\lambda }}{(\partial _{\mathbf{n}}w)f\ dv_{g}^{\prime }}],
\end{equation}%
where $\partial _{\mathbf{n}}$ is the unit outward normal derivative with
respect to $g$, and $v_{g}^{\prime }$ is the induced volume form on $%
S_{\lambda }$.
\end{lemma}

\noindent \emph{Proof.} This is proved through direct computation. Because $%
f $ is constant on $S_{\lambda }$,
\begin{equation}
\Delta _{g}f=\partial _{n}^{2}f+H\partial _{\mathbf{n}}f  \label{hun}
\end{equation}%
on $S_{\lambda }$, where $H$ is the mean curvature of $S_{\lambda }\subset
U_{\lambda }$ with respect to $g$. Also notice that on $S_{\lambda }$,%
\begin{equation}
{\frac{d}{{d\lambda }}}dv_{g}^{\prime }=Hdv_{g}^{\prime }.  \label{gry}
\end{equation}%
Thus, if $S_{\lambda }$ is smooth, which is true for almost all $\lambda $,
\begin{eqnarray}
{\frac{d}{{d\lambda }}}\int\limits_{S_{\lambda }}{(\partial _{\mathbf{n}%
}\lambda )f\ dv_{g}^{\prime }} &=&\int\limits_{S_{\lambda }}{{\frac{{%
(\partial _{n}^{2}\lambda )f+(\partial _{\mathbf{n}}\lambda )(\partial _{%
\mathbf{n}}f)+(\partial _{\mathbf{n}}\lambda )fH}}{{\partial _{\mathbf{n}%
}\lambda }}}\ dv_{g}^{\prime }}  \notag \\
&=&\int\limits_{S_{\lambda }}{({\frac{{\Delta _{g}\lambda }}{{\partial _{%
\mathbf{n}}\lambda }}}f+\partial _{\mathbf{n}}f)\ dv_{g}^{\prime }},
\end{eqnarray}%
by (\ref{hun}) and (\ref{gry}). Using Stokes' Theorem and the co-area
formula, we get
\begin{eqnarray}
\int\limits_{U_{\lambda }}{\Delta _{g}f\ dv_{g}} &=&-\int\limits_{S_{\lambda
}}{(\partial _{\mathbf{n}}f)\ dv_{g}^{\prime }}=-{\frac{d}{{d\lambda }}}%
\int\limits_{S_{\lambda }}{(\partial _{\mathbf{n}}\lambda )f\ dv_{g}^{\prime
}}+\int\limits_{S_{\lambda }}{{\frac{{\Delta \lambda }}{{\partial _{\mathbf{n%
}}\lambda }}}fdv_{g}^{\prime }}  \notag \\
&=&-\lambda {\frac{d}{{d\lambda }}}\int\limits_{S_{\lambda }}{({\frac{{%
\partial _{\mathbf{n}}\lambda }}{\lambda }})f\ dv_{g}^{\prime }}%
-\int\limits_{S_{\lambda }}{({\frac{{\partial _{\mathbf{n}}\lambda }}{%
\lambda }})fdv_{g}^{\prime }}+\lambda {\frac{d}{{d\lambda }}}%
\int\limits_{U_{\lambda }}{({\frac{{\Delta _{g}\lambda }}{\lambda }})f\
dv_{g}}.
\end{eqnarray}%
The co-area formula also leads to
\begin{equation}
\int\limits_{S_{\lambda }}{({\frac{{\partial _{\mathbf{n}}\lambda }}{\lambda
}})f\ dv_{g}^{\prime }}=\lambda {\frac{d}{{d\lambda }}}\int\limits_{U_{%
\lambda }}{({\frac{{(\partial _{\mathbf{n}}\lambda )^{2}}}{{\lambda ^{2}}}}%
)f\ dv_{g}}.
\end{equation}%
Finally, we have
\begin{equation*}
\int\limits_{U_{\lambda }}{\Delta _{g}f\ dv_{g}}=\lambda {\frac{d}{{d\lambda
}}}[\int\limits_{U_{\lambda }}{({\frac{\Delta _{g}\lambda }{\lambda }}-{%
\frac{(\partial _{\mathbf{n}}\lambda )^{2}}{\lambda ^{2}}})\ dv_{g}}%
-\int\limits_{S_{\lambda }}{{\frac{\partial _{\mathbf{n}}\lambda }{\lambda }}%
f\ dv_{g}^{\prime }}]
\end{equation*}%
\begin{equation*}
=\lambda {\frac{d}{{d\lambda }}}[\int\limits_{U_{\lambda }}{(\Delta
_{g}w)fdv_{g}}-\int\limits_{S_{\lambda }}{(\partial _{\mathbf{n}}w)f\
dv_{g}^{\prime }}]
\end{equation*}%
by the fact that $\lambda =e^{w}$ on $S_{\lambda }$.\qed

We are now in the position to prove Theorem~\ref{4-1}. This is an
interesting generalization of the main results of [CQY2], where we explore
heavily the (local) divergence structure of the Pfaffian curvature..

{\noindent \emph{Proof of Theorem~\ref{4-1}.} } As in the proof of Lemma~\ref%
{pf-unique}, we first apply Lemma~\ref{savior} to show
\begin{equation}
\|\nabla w\|\leq C e^w.
\end{equation}

For notational simplicity, we use upper-bar to to denote geometric objects
constructed with $g_{w}$ metric.

Notice by (\ref{faf1}) and (\ref{faf2}),
\begin{equation}
\sigma _{2}=e^{-4w}\Delta ^{2}w+\bar{\Delta}J  \label{suan0}
\end{equation}%
Apply Lemma~\ref{reduce} with the metric $g_{0}$, we have
\begin{equation}
\int\limits_{U_{\lambda }}{\Delta ^{2}wdv}=\lambda \frac{d}{{d\lambda }}%
[\int\limits_{U_{\lambda }}{(\Delta w)^{2}dv}-\int\limits_{S_{\lambda }}{%
(\partial _{\mathbf{n}}w)\Delta wdv^{\prime }}]  \label{suan1}
\end{equation}%
It's a simple observation that Lemma~\ref{reduce} can also be applied with
respect to the metric $g=g_{w}$ to get
\begin{equation}
\int\limits_{U_{\lambda }}{\bar{\Delta}Jdv_{g}}=\lambda \frac{d}{{d\lambda }}%
[\int\limits_{U_{\lambda }}{\bar{\Delta}wJdv_{g}}-\int\limits_{S_{\lambda }}{%
(\bar{\nabla}_{\bar{n}}w)Jdv_{g}^{\prime }}],  \label{suan2}
\end{equation}%
with $\bar{\mathbf{n}}=e^{-w}\mathbf{n}$ being the unit normal vector of $%
S_{\lambda }=\partial U_{\lambda }$ with respect to $g$. Notice that $%
J=-e^{-2w}(\Delta w+\Vert \nabla w\Vert ^{2})$ and $\bar{\nabla}_{\mathbf{%
\bar{n}}}w=e^{-w}\nabla _{\mathbf{n}}w$, we get, from (\ref{suan0}), (\ref%
{suan1}) and (\ref{suan2}), 

\begin{eqnarray}
\int\limits_{U_{\lambda }}{\sigma _{2}dv_{g}} &=&\lambda \frac{d}{{d\lambda }%
}\int\limits_{U_{\lambda }}{[(\Delta w)^{2}-(\Delta w+2\Vert \nabla w\Vert
^{2})(\Delta w+\Vert \nabla w\Vert ^{2})]dv}  \notag \\
&&-\lambda \frac{d}{{d\lambda }}\int\limits_{S_{\lambda }}{\partial _{%
\mathbf{n}}w[\Delta w-(\Delta w+\Vert \nabla w\Vert ^{2})]dv^{\prime }}.
\notag \\
&=&\lambda \frac{d}{{d\lambda }}[\int\limits_{U_{\lambda }}{(3Je^{2w}+\Vert
\nabla w\Vert ^{2})\Vert \nabla w\Vert ^{2}dv}+\int\limits_{S_{\lambda }}{%
(\partial _{\mathbf{n}}w)\Vert \nabla w\Vert ^{2}dv^{\prime }}].
\label{5.28'}
\end{eqnarray}%
Define
\begin{equation*}
F(\lambda )\equiv \int\limits_{U_{\lambda }}{(3Je^{2w}+\Vert \nabla w\Vert
^{2})\Vert \nabla w\Vert ^{2}dv}+\int\limits_{S_{\lambda }}{(\partial _{%
\mathbf{n}}w)\Vert \nabla w\Vert ^{2}dv^{\prime }}.
\end{equation*}%
Notice $J$ and $\partial _{\mathbf{n}}w$ both being non-negative,
\begin{equation*}
F(\lambda )\geq C\int\limits_{U_{\lambda }}{\Vert \nabla w\Vert ^{2}e^{2w}dv}%
=C\int\limits_{U_{\lambda }}{\Vert \nabla e^{w}\Vert ^{2}dv}%
=C[\int\limits_{U_{\lambda }}{(-\Delta \lambda )\cdot \lambda dv}%
+\int\limits_{S_{\lambda }}{\partial _{\mathbf{n}}\lambda \cdot \lambda
dv^{\prime }}]
\end{equation*}%
\begin{equation}
=C[\int\limits_{U_{\lambda }}{J\lambda ^{4}dv}+\int\limits_{S_{\lambda }}{%
\partial _{\mathbf{n}}\lambda \cdot \lambda dv^{\prime }}]\geq
C\int\limits_{U_{\lambda }}{e^{4w}dv}.
\end{equation}

Following [CQY2], we now apply Lemma~\ref{blow-up} and Lemma~2.6 of [CQY2]
to get
\begin{eqnarray}
F(\lambda ) &\leq &-C(\lambda ^{{\frac{{3}}{4}}\dim {\Lambda }}-1),\ \
\mathrm{{if}\ \dim \Lambda >0;} \\
F(\lambda ) &\leq &-N\ln \lambda ,\ \ \mathrm{{if}\ \dim \Lambda =0,\
H^{0}(\Lambda )=\infty ,}
\end{eqnarray}%
for any large integer $N.$ Hence, there is a sequence of $\lambda _{i}$,
such that $\lambda _{i}\rightarrow \infty $ and $\lambda _{i}{\frac{d}{{%
d\lambda }}}F(\lambda _{i})\rightarrow -\infty $ as $i$ tends to infinity.
But by (\ref{5.28'}) this contradicts with assumption (A4). Hence we have
proved that $S^{4}\backslash \Omega $ is of Hausdorff dimension 0 and it is
actually finite. The proof is finished.

\qed


\end{document}